\documentclass[11pt]{article}

\usepackage[T1]{fontenc}
\usepackage[margin=1in]{geometry}
\usepackage{amsmath,amssymb,amsthm}
\usepackage{booktabs}
\usepackage{array}
\usepackage{enumitem}
\usepackage{xcolor}
\usepackage{hyperref}

\hypersetup{
  colorlinks=true,
  linkcolor=blue!55!black,
  citecolor=blue!55!black,
  urlcolor=blue!55!black,
  pdftitle={Non-Pseudo-Arithmetic Hyperbolic 4-Manifolds},
  pdfauthor={Alexander Kolpakov}
}

\numberwithin{equation}{section}

\newtheorem{theorem}{Theorem}[section]
\newtheorem{proposition}[theorem]{Proposition}
\newtheorem{lemma}[theorem]{Lemma}
\newtheorem{corollary}[theorem]{Corollary}

\theoremstyle{definition}
\newtheorem{definition}[theorem]{Definition}

\theoremstyle{remark}
\newtheorem{remark}[theorem]{Remark}

\newcommand{\Q}{\mathbb{Q}}
\newcommand{\R}{\mathbb{R}}
\newcommand{\C}{\mathbb{C}}
\newcommand{\HH}{\mathbb{H}}

\newcommand{\tr}{\operatorname{tr}}
\newcommand{\Ad}{\operatorname{Ad}}
\newcommand{\Isom}{\operatorname{Isom}}
\newcommand{\PO}{\operatorname{PO}}

\newcommand{\minpoly}{\operatorname{minpoly}}

\newcommand{\vol}{\operatorname{vol}}
\newcommand{\code}[1]{\texttt{#1}}

\title{Non-Pseudo-Arithmetic Hyperbolic \(4\)-Manifolds}
\author{Alexander Kolpakov}
\date{30 October 2026}

\begin{document}

\maketitle

\begin{abstract}
We give three explicit finite-volume hyperbolic Coxeter lattices in dimension
four whose adjoint trace fields are not totally real.  They are the reflection
groups of the Ma--Zheng polytopes
\(P_{11,6},P_{11,7},P_{11,8}\).  The obstruction is contained in two cyclic
products of the doubled unit-normal Gram matrix \(A=2G\).  Those products
recover elements \(s=\sqrt5\) and \(u\) satisfying
\[
 (s-1)u^2+2(s-3)u+s-7=0.
\]
After the embedding \(s\mapsto-\sqrt5\), this quadratic has discriminant
\(8(1-\sqrt5)<0\).  The cyclic-product field therefore has a complex place.
Vinberg's identification of that field with the adjoint trace field and the
Emery--Mila total-reality obstruction imply that the three lattices are not
pseudo-arithmetic.  Torsion-free orientation-preserving finite-index
subgroups yield orientable noncompact finite-volume hyperbolic
\(4\)-manifolds with the same obstruction.

Exact computations check the census data and provide independent field and
finite-volume certificates.  A pinned Lean~4 development verifies the
radical identity and the number-field obstruction.  It also checks the
lattice and manifold implications conditional on explicit geometric and
literature hypotheses.
\end{abstract}

\tableofcontents

\section{Introduction}
\label{sec:introduction}

\subsection{The question and the answer}

A \emph{lattice} in the isometry group of real hyperbolic \(n\)-space is a
discrete group whose quotient has finite hyperbolic volume.  Arithmetic
lattices are constructed from algebraic groups over number fields.  Many
important nonarithmetic lattices are assembled by gluing arithmetic pieces.
Emery and Mila introduced \emph{pseudo-arithmeticity} to give one common
algebraic framework for arithmetic lattices, quasi-arithmetic lattices, and
the familiar gluing constructions \cite{EmeryMila}.  They asked whether every
lattice in \(\PO(n,1)\), for \(n>3\), is pseudo-arithmetic and, more weakly,
whether every such lattice has totally real adjoint trace field
\cite[Questions~1.8 and~1.10]{EmeryMila}.

Both questions have negative answers in dimension four.

\begin{theorem}[Main theorem]
\label{thm:main}
For \(j\in\{6,7,8\}\), let \(P_{11,j}\) be the finite-volume hyperbolic
Coxeter \(4\)-polytope in the Ma--Zheng census recorded in
Table~\ref{tab:three-polytopes}, and let
\(\Gamma_{11,j}<\PO(4,1)\) be its reflection group.  Then the adjoint trace
field of \(\Gamma_{11,j}\) is not totally real.  Consequently
\(\Gamma_{11,j}\) is not pseudo-arithmetic.
\end{theorem}

\begin{corollary}[Manifold form]
\label{cor:manifolds}
There exist orientable, noncompact, finite-volume hyperbolic
\(4\)-manifolds whose fundamental groups are not pseudo-arithmetic.
\end{corollary}

Ma and Zheng classify all \(331\) finite-volume hyperbolic Coxeter
\(4\)-polytopes with seven facets \cite{MaZheng}.  The present examples are
three consecutive entries in that census.  The proof needs only four entries
of their Gram matrices and one common radical.  It does not require a
numerical approximation, a primitive element for the entire field, or a
classification of all embeddings of that field.

Emery--Mila checked the higher-dimensional Coxeter examples then available
to them.  The final 2026 version of
Belolipetsky--Bogachev--Kolpakov--Slavich still described the known examples
in \(\PO(n,1)\), \(n>3\), as pseudo-arithmetic \cite{BBKS}.  A targeted
literature search through 24 September 2026 found no earlier recorded
non-pseudo-arithmetic lattice in these dimensions.  We know of no previous
example.

\subsection{Definitions and conventions}
\label{subsec:primer}

We recall the definitions used in the argument.

\paragraph{Hyperbolic space and lattices.}
On \(\R^{n+1}\), put
\[
 \langle x,y\rangle_{n,1}
   =x_1y_1+\cdots+x_ny_n-x_{n+1}y_{n+1}.
\]
The upper sheet
\[
 \HH^n=\{x\in\R^{n+1}:\langle x,x\rangle_{n,1}=-1,
                       \ x_{n+1}>0\}
\]
is the hyperboloid model of real hyperbolic space.  The subgroup of the
orthogonal group preserving this sheet acts by isometries and is naturally
isomorphic to \(\PO(n,1)=\mathrm O(n,1)/\{\pm I\}\).  A subgroup
\(\Gamma<\Isom(\HH^n)\) is \emph{discrete} if the identity is isolated in
\(\Gamma\).  A discrete subgroup is a \emph{lattice} if
\(\Gamma\backslash\HH^n\) has finite volume.  The quotient is an orbifold in
general and a manifold when \(\Gamma\) is torsion-free.

\paragraph{Coxeter polytopes and Gram matrices.}
A convex hyperbolic polytope is \emph{Coxeter} when every nonzero dihedral
angle is \(\pi/m\) for an integer \(m\ge2\); zero angles at ideal intersections
are also allowed.  Reflections in its facets generate a
discrete group for which the polytope is a fundamental domain.  If the
polytope has finite volume, this reflection group is a lattice.  Choose an
outward unit spacelike normal \(e_i\) to each facet.  The matrix
\(G=(\langle e_i,e_j\rangle_{n,1})\) is the \emph{unit-normal Gram matrix}.
Its diagonal entries are one.  Off the diagonal one has
\[
 g_{ij}=
 \begin{cases}
 -\cos(\pi/m_{ij}),&\text{the facets meet at angle }\pi/m_{ij},\\
 -1,&\text{the facets are parallel},\\
 -\cosh(\rho_{ij}),&\text{they are disjoint at distance }\rho_{ij}.
 \end{cases}
\]
We use \(A=2G\).  A \emph{cyclic product} is a product
\(A_{i_0i_1}A_{i_1i_2}\cdots A_{i_mi_0}\).  The field generated over \(\Q\)
by every such product is the \emph{Vinberg field} \(k(P)\).

\paragraph{Number fields and total reality.}
A \emph{number field} is a finite extension of \(\Q\).  An embedding of a
number field \(K\) is an injective field homomorphism \(K\hookrightarrow\C\)
fixing \(\Q\).  The field \(K\) is \emph{totally real} when every such
embedding has image in \(\R\).  Thus one nonreal embedding is a complete
certificate that \(K\) is not totally real.  We use the standard extension
theorem: if \(L\subset K\) are number fields, every embedding
\(L\hookrightarrow\C\) extends to an embedding \(K\hookrightarrow\C\)
(Proposition~\ref{prop:embedding-extension}).

\paragraph{Adjoint trace fields.}
For a matrix Lie group, conjugation differentiates to a linear action
\(\Ad\) on its Lie algebra.  The \emph{adjoint trace field} of a lattice
\(\Gamma\) is
\[
 k(\Gamma)=
 \Q\bigl(\{\tr(\Ad\gamma):\gamma\in\Gamma\}\bigr)\subset\R.
\]
It is unchanged on passing to a finite-index subgroup.  Vinberg's
field-of-definition theorem identifies this field, for a finite-volume
Coxeter reflection group, with the cyclic-product field \(k(P)\).  This
identification is an imported theorem, not a new result proved here.

\paragraph{Pseudo-arithmeticity.}
Vinberg associates to a Zariski-dense lattice \(\Gamma\) an algebraic group
over \(k(\Gamma)\), called its ambient group.  An algebraic group over a
totally real field \(k\) is \emph{admissible} for hyperbolic \(n\)-space when
one real factor is the required noncompact projective orthogonal group and
the other real factors are compact.  In Emery--Mila's definition, a lattice
is \emph{pseudo-arithmetic} when its ambient group is obtained from such an
admissible group by scalar extension through a finite totally real multiquadratic
extension \(K/k\).  Their definition implies that \(K=k(\Gamma)\) is totally
real.  This last necessary condition is the only part of
pseudo-arithmeticity used in the proof.

\paragraph{Formal verification terminology.}
Lean checks a theorem only from its declarations, imported library theorems,
and explicitly supplied hypotheses.  In this release, the new algebraic
argument is proved inside Lean.  A literature result such as Vinberg's theorem
is represented by a named hypothesis with the same logical type needed by
the manuscript.  Lean checks every implication after that boundary, but this
does not constitute a new formal proof of Vinberg's theorem.  No
project-defined \code{axiom}, \code{opaque}, \code{sorry}, or \code{admit}
is used.

\subsection{The certificate in one paragraph}
\label{subsec:certificate-preview}

Let \(G\) be any of the three unit-normal Gram matrices and put \(A=2G\).
With zero-based facet indices, the common entries
\[
 A_{12}=-\frac{1+\sqrt5}{2},\qquad
 A_{03}=-1,\qquad A_{35}=-\sqrt2,\qquad A_{50}=-2a
\]
give two cyclic products
\[
 c_{12}=A_{12}A_{21}=\frac{3+\sqrt5}{2},\qquad
 c_{035}=A_{03}A_{35}A_{50}=-2\sqrt2\,a.
\]
Hence the Vinberg field contains
\(s=2c_{12}-3=\sqrt5\) and
\(u=-c_{035}/2=\sqrt2a\).  The exact radical for \(a\) implies
\[
 (s-1)u^2+2(s-3)u+s-7=0.
\]
At the conjugate embedding of \(\Q(\sqrt5)\), this quadratic has negative
discriminant.  It therefore produces a nonreal embedding of
\(\Q(s,u)\), which extends to the whole Vinberg field.  This is the entire
new obstruction.

\subsection{Logical status of the ingredients}
\label{subsec:logical-status}

We state the ingredients of the proof and their logical dependencies
explicitly for the subsequent Lean~4 formalization.  We distinguish four
layers.

\begin{enumerate}[label=\textbf{L\arabic*.},leftmargin=3em]
\item \emph{Elementary proof:} the matrix entries, two products, radical
  elimination, discriminant, and complex-place deduction are written out in
  Sections~\ref{sec:coxeter-data}--\ref{sec:main-proof}.
\item \emph{Exact computation:} independent Python and PARI/GP programs check
  the source transcription, rank, inertia, field signatures, and finite-volume
  combinatorics.  They support the input data but are not needed to trust the
  two-cycle algebra once the displayed data are accepted.
\item \emph{Lean kernel checking:} the internal implication from the radical
  formula and cycle certificate to non-total reality and non-pseudo-
  arithmeticity is formalized.  The finite-index corollary is checked from
  explicitly named geometric hypotheses.
\item \emph{External mathematics:} the census classification, Vinberg's
  field identification, the Emery--Mila necessary condition, Selberg's lemma,
  and finite-index invariance are cited rather than reproved.  The exact list
  appears in Section~\ref{sec:external-register}.
\end{enumerate}

This breakdown identifies the statements checked in Lean and the external
results supplied as hypotheses.

\section{Preliminaries}
\label{sec:basics}

\subsection{Embeddings and complex places}

We state the following standard field-theoretic facts explicitly for the
subsequent Lean~4 formalization.  We write \(\minpoly_\Q(\alpha)\) for the
minimal polynomial of an algebraic number \(\alpha\) over \(\Q\).

\begin{lemma}[Quadratic conjugation]
\label{lem:sqrt5-minpoly}
One has
\[
 \minpoly_\Q(\sqrt5)=X^2-5.
\]
The two \(\Q\)-embeddings of \(\Q(\sqrt5)\) into \(\C\) are
determined by \(\sqrt5\mapsto\pm\sqrt5\).
\end{lemma}

\begin{proposition}[Extension of embeddings]
\label{prop:embedding-extension}
Let \(L\subset K\) be number fields and let
\(\sigma:L\hookrightarrow\C\) be an embedding.  Then \(\sigma\) extends to
an embedding \(\widetilde\sigma:K\hookrightarrow\C\).
\end{proposition}

\begin{proof}
Because number fields have characteristic zero, \(K/L\) is finite and
separable.  By the primitive-element theorem there is \(\theta\in K\) with
\(K=L(\theta)\).  Let \(p\in L[X]\) be the minimal polynomial of \(\theta\).
Apply \(\sigma\) coefficientwise to obtain \(p^\sigma\in\C[X]\).  The
fundamental theorem of algebra supplies a root \(z\in\C\) of
\(p^\sigma\).  The map
\[
 L[X]\longrightarrow\C,
 \qquad \sum_i a_iX^i\longmapsto\sum_i\sigma(a_i)z^i
\]
has \((p)\) in its kernel and therefore descends to a homomorphism
\(L[X]/(p)\to\C\).  Identifying \(L[X]/(p)\) with \(L(\theta)=K\) gives the
extension.  A unital homomorphism from a field has zero kernel, so it is an
embedding.
\end{proof}

\begin{lemma}[Negative discriminant]
\label{lem:quadratic-test}
A quadratic polynomial \(aX^2+bX+c\in\R[X]\), with \(a\ne0\) and
\(b^2-4ac<0\), has no real root.
\end{lemma}

\begin{lemma}[Persistence of a complex place]
\label{prop:subfield-persists}
For number fields \(L\subset K\), every nonreal embedding of \(L\)
extends to a nonreal embedding of \(K\).  In particular, if \(L\) is not
totally real, then neither is \(K\).
\end{lemma}

\subsection{The hyperboloid model and reflecting hyperplanes}

Let \(V=\R^{n+1}\) with the form \(\langle\ ,\ \rangle_{n,1}\) defined in
Subsection~\ref{subsec:primer}.  A nonzero vector is spacelike, null, or
timelike according as its squared Lorentzian norm is positive, zero, or
negative.  If \(e\) is a unit spacelike vector, then
\[
 H_e=\{x\in\HH^n:\langle x,e\rangle_{n,1}=0\}
\]
is a totally geodesic hyperplane, and
\[
 r_e(x)=x-2\langle x,e\rangle_{n,1}e
\]
is the reflection fixing \(H_e\) pointwise.  Indeed,
\(r_e(e)=-e\), \(r_e\) fixes the orthogonal complement of \(e\), and direct
expansion shows
\(\langle r_e(x),r_e(y)\rangle_{n,1}=\langle x,y\rangle_{n,1}\).

For a polytope described as the intersection of half-spaces
\(\langle x,e_i\rangle_{n,1}\le0\), the vectors \(e_i\) are outward normals.
The inner products of pairs of normals give the Gram formulas in
Subsection~\ref{subsec:primer}; these are the standard two-dimensional
Lorentzian calculations for intersecting, asymptotic, and ultraparallel
hyperplanes \cite[Chapter~3]{Ratcliffe}.

\begin{definition}
A Coxeter matrix is symmetric, has diagonal entries one, and has off-diagonal labels
\(m_{ij}\in\{2,3,\ldots\}\cup\{\infty\}\).  Its Coxeter group is generated
by involutions \(r_i\) with relations \((r_ir_j)^{m_{ij}}=1\) when
\(m_{ij}<\infty\).  For a hyperbolic Coxeter polytope, these generators are
the geometric facet reflections.
\end{definition}

The fact that the geometric reflection group is discrete and that the
polytope is a fundamental domain is the Poincar\'e--Coxeter theorem.  We use
it as standard background.  Finite volume of the fundamental domain is
exactly finite covolume of the group.

\subsection{From an orbifold to a manifold}
\label{subsec:orbifold-manifold}

Two elementary finite-index observations will be used.  First, the kernel of
the orientation character \(\Gamma\to\{\pm1\}\) has index at most two.
For a group generated by hyperplane reflections, each generator reverses
orientation, so the kernel has index exactly two.  Second, a subgroup of
finite index \(d\) gives a degree-\(d\) orbifold cover and multiplies volume
by \(d\).  Thus finite volume is preserved in both directions.

The following external theorem removes torsion.

\begin{theorem}[Selberg's lemma]
\label{thm:selberg}
Every finitely generated subgroup of \(\mathrm{GL}_N(F)\), where \(F\) has
characteristic zero, has a torsion-free subgroup of finite index.
\end{theorem}

We cite Selberg's original result \cite{Selberg}.  A Coxeter group with
finitely many facets is finitely generated, and its reflection
representation is linear over a characteristic-zero field, so the theorem
applies.  A torsion-free discrete group acts freely on \(\HH^n\): an element
fixing a point lies in a compact point stabilizer, and discreteness makes the
cyclic subgroup finite, hence trivial.  The quotient is therefore a
manifold.  Passing inside the orientation-preserving subgroup makes it
orientable.  Finally, a finite cover of a noncompact finite-volume quotient
is noncompact: if the cover were compact, its continuous image would be
compact, contrary to noncompactness of the base.

\section{Adjoint trace fields and pseudo-arithmeticity}
\label{sec:trace-fields}

\subsection{The adjoint representation and ambient group}

Let \(\mathbf G\) be a real algebraic matrix group with Lie algebra
\(\mathfrak g\).  Conjugation by \(g\in\mathbf G(\R)\) maps
\(X\in\mathfrak g\) to \(gXg^{-1}\).  This linear map is \(\Ad(g)\).  Its
trace is independent of a choice of basis for \(\mathfrak g\).  For a
Zariski-dense subgroup \(\Gamma\), define
\[
 k(\Gamma)=\Q\bigl(\tr(\Ad\gamma):\gamma\in\Gamma\bigr).
\]
For lattices in a connected noncompact simple Lie group, Borel density gives
the required Zariski density.  Vinberg's theory produces a form of the
ambient adjoint algebraic group over \(k(\Gamma)\), unique up to
\(k(\Gamma)\)-isomorphism \cite{VinbergRings}.  We do not need its
construction; only the field is used below.

We also use the standard finite-index invariance
\begin{equation}
 [\Gamma:\Lambda]<\infty
 \quad\Longrightarrow\quad k(\Lambda)=k(\Gamma).
\label{eq:finite-index-field}
\end{equation}
This is part of Vinberg's commensurability theory.  One inclusion is
immediate from \(\Lambda\subset\Gamma\); equality is the nontrivial external
statement and is listed as Input~E6 in
Section~\ref{sec:external-register}.

\subsection{Pseudo-admissibility}

An algebraic \(k\)-group \(\mathbf H\) is \emph{admissible for
\(\PO(n,1)\)} when \(k\) is totally real and, among the groups obtained from
the real embeddings of \(k\), exactly one has the required noncompact
factor \(\PO(n,1)^\circ\), while the others are compact.  This is the usual
arithmetic admissibility condition.

\begin{definition}[Emery--Mila]
Let \(n>3\).  A group over a number field \(K\) is
\emph{pseudo-admissible over \(K/k\)} if it is the scalar extension
\(\mathbf H_K\) of an admissible \(k\)-group through a finite, possibly
trivial, totally real multiquadratic extension
\[
 K=k(\sqrt{d_1},\ldots,\sqrt{d_t}).
\]
A lattice \(\Gamma<\PO(n,1)\) is \emph{pseudo-arithmetic} when its ambient
group is pseudo-admissible \cite[Definitions~1.6--1.7]{EmeryMila}.
\end{definition}

The trivial case \(t=0\) includes quasi-arithmetic lattices.
Total reality means that every embedding of \(K\) is real, as required in
\cite[Section~1.4]{EmeryMila}.
Equivalently, the generators \(d_i\in k^\times\) can be chosen totally
positive, meaning \(\sigma(d_i)>0\) for every embedding
\(\sigma:k\hookrightarrow\R\).  Merely requiring \(K\subset\R\) at one
embedding would not imply total reality.

\begin{proposition}[Total-reality obstruction]
\label{prop:total-real-obstruction}
If \(n>3\) and \(\Gamma<\PO(n,1)\) is pseudo-arithmetic, then its adjoint
trace field \(k(\Gamma)\) is totally real.
\end{proposition}

\begin{proof}
We use the necessary condition stated by Emery--Mila immediately before
Question~1.10 in \cite[Section~1.4]{EmeryMila}.  In the notation above, an
embedding \(\tau:K\hookrightarrow\C\) restricts to a real embedding
\(\sigma\) of \(k\), and
\(\tau(\sqrt{d_i})^2=\sigma(d_i)>0\).  Hence every generator has real
image, so \(\tau(K)\subset\R\).  Since \(K=k(\Gamma)\), the conclusion
follows.
\end{proof}

The contrapositive is what we use:
\begin{equation}
 k(\Gamma)\text{ not totally real}
 \quad\Longrightarrow\quad
 \Gamma\text{ not pseudo-arithmetic}.
\label{eq:total-real-contrapositive}
\end{equation}

\subsection{Vinberg's cyclic-product field}
\label{subsec:vinberg}

Let \(P\subset\HH^n\) be a finite-volume Coxeter polytope with unit-normal
Gram matrix \(G=(g_{ij})\), and put \(A=2G\).  For a closed sequence of facet
indices \(i_0,i_1,\ldots,i_m,i_0\), its cyclic product is
\[
 A_{i_0i_1}A_{i_1i_2}\cdots A_{i_{m-1}i_m}A_{i_mi_0}.
\]
Repeated indices are allowed, and reversing the cycle gives the same value
because \(A\) is symmetric.  Define
\[
 k(P)=\Q(\text{all cyclic products in }A).
\]

\begin{theorem}[Vinberg]
\label{thm:vinberg-cycles}
If \(P\) has finite volume and \(\Gamma(P)\) is its reflection group, then
\[
 k(\Gamma(P))=k(P).
\]
\end{theorem}

This is the cyclic-product form of Vinberg's field-of-definition theorem
\cite{VinbergRings}; see also \cite[Section~4.1]{EmeryMila} and
\cite[Section~3.1]{GuoMaZangZheng}.  It is an external input.  The factor two
fixes the normalization of the displayed cyclic products: starting from unit
normals, Vinberg's diagonal-two convention gives \(A=2G\).  The generated
field is unchanged if one uses \(G\), since a cycle of length \(\ell\) is
multiplied by the rational factor \(2^\ell\).

\begin{remark}[Entry field versus trace field]
The field generated by all individual entries of \(G\) can be strictly
larger than \(k(P)\).  The relevant rescaling invariance concerns the
reflection coefficients
\(a_{ij}=2\langle e_i,e_j\rangle/\langle e_i,e_i\rangle\), which equal
\(A_{ij}\) for unit normals.  Replacing \(e_i\) by \(\lambda_i e_i\)
changes \(a_{ij}\) to \((\lambda_j/\lambda_i)a_{ij}\), so the factors
cancel around a closed cycle.  Raw bilinear Gram entries do not have this
rescaling rule.  Thus a radical in an individual Gram entry need not lie in
the adjoint trace field; the field inclusions below follow from explicit
cyclic products.
\end{remark}

\section{The three Coxeter polytopes and their exact matrices}
\label{sec:coxeter-data}

\subsection{How the census vector is read}

We number the seven facets \(0,1,\ldots,6\).  The \(21\) unordered pairs are
listed in lexicographic blocks:
\begin{equation}
\begin{split}
 &(01,02,03,04,05,06)\mid(12,13,14,15,16)\mid(23,24,25,26)\\
 &\hspace{43mm}\mid(34,35,36)\mid(45,46)\mid(56).
\end{split}
\label{eq:pair-order}
\end{equation}
An integer \(m\ge2\) means dihedral angle \(\pi/m\), the symbol \(0\)
means parallel facets, and a letter means an ultraparallel pair.  Thus the
initial block \(2232ab\) means
\[
 m_{01}=2,\quad m_{02}=2,\quad m_{03}=3,\quad m_{04}=2,
 \quad g_{05}=-a,\quad g_{06}=-b.
\]

The relevant published rows are as follows.  The displayed volume is the
volume of the Coxeter fundamental polytope.

\begin{table}[ht]
\centering
\small
\begin{tabular}{@{}c l c c@{}}
\toprule
polytope & Coxeter vector & cusps & volume \\
\midrule
\(P_{11,6}\) & \(\mathtt{2232ab\ 52322\ 2220\ 442\ 22\ 2}\)
    & \(1\) & \(13\pi^2/864\) \\
\(P_{11,7}\) & \(\mathtt{2232ab\ 52323\ 2220\ 442\ 22\ 2}\)
    & \(1\) & \(\pi^2/54\) \\
\(P_{11,8}\) & \(\mathtt{2232ab\ 52324\ 2220\ 442\ 22\ 2}\)
    & \(2\) & \(19\pi^2/864\) \\
\bottomrule
\end{tabular}
\caption{The three Ma--Zheng census rows used in the proof.}
\label{tab:three-polytopes}
\end{table}

The assertion that these rows define finite-volume polytopes is external
Input~E1.  Section~\ref{sec:certification} describes independent exact checks,
including a direct finite-volume check for \(P_{11,8}\).

\subsection{Radicals and a uniform matrix}

Put
\[
 s=\sqrt5,\qquad r=\sqrt2,
\]
and let \(q_6=2,q_7=3,q_8=4\).  The common first dotted entry is
\begin{equation}
 a=\frac12\sqrt{7+s+2\sqrt{2+2s}}.
\label{eq:a-radical}
\end{equation}
The second dotted entries are
\begin{equation}
\begin{aligned}
 b_6&=\sqrt{\frac92+2s+\sqrt{38+17s}},\\
 b_7&=r+\frac{\sqrt{10}}2+\frac12\sqrt{29+13s},\\
 b_8&=\frac14\left(2+3r+2s+rs
       +2\sqrt{19+14r+9s+6rs}\right).
\end{aligned}
\label{eq:b-radicals}
\end{equation}
All square roots mean the positive real root.  Their radicands are positive,
and the resulting \(a,b_j\) exceed one, as required for ultraparallel
facets.  Only \(a\) is used in the obstruction.

For reference, the complete unit-normal Gram matrix is
\begin{equation}
G_j=
\begin{pmatrix}
1&0&0&-\tfrac12&0&-a&-b_j\\
0&1&-\tfrac{1+s}{4}&0&-\tfrac12&0&-\cos(\pi/q_j)\\
0&-\tfrac{1+s}{4}&1&0&0&0&-1\\
-\tfrac12&0&0&1&-\tfrac r2&-\tfrac r2&0\\
0&-\tfrac12&0&-\tfrac r2&1&0&0\\
-a&0&0&-\tfrac r2&0&1&0\\
-b_j&-\cos(\pi/q_j)&-1&0&0&0&1
\end{pmatrix}.
\label{eq:full-gram}
\end{equation}
To check the transcription, read each row against \eqref{eq:pair-order}.
For example, \(m_{12}=5\) gives
\(g_{12}=-\cos(\pi/5)=-(1+s)/4\); \(m_{26}=0\) gives \(g_{26}=-1\);
and \(m_{34}=m_{35}=4\) gives \(-r/2\).  The elementary special-angle
identities used here are
\[
 \cos(\pi/3)=\frac12,\qquad
 \cos(\pi/4)=\frac{\sqrt2}{2},\qquad
 \cos(\pi/5)=\frac{1+\sqrt5}{4}.
\]
The first two follow from an equilateral triangle and an isosceles right
triangle.  For the third, put \(z=e^{i\pi/5}\) and
\(x=z+z^{-1}=2\cos(\pi/5)\).  Since \(z^5=-1\) and \(z\ne-1\), division of
\(z^5+1\) by \(z+1\) gives
\(z^4-z^3+z^2-z+1=0\).  After division by \(z^2\), this becomes
\[
 (z^2+z^{-2})-(z+z^{-1})+1=0.
\]
Because \(z^2+z^{-2}=x^2-2\), we get \(x^2-x-1=0\).  Positivity selects
\(x=(1+\sqrt5)/2\).

\subsection{The two cyclic products}
\label{subsec:two-cycles}

Put \(A_j=2G_j\).  Four common entries of \(A_j\) are
\begin{equation}
 A_{12}=-\frac{1+s}{2},\qquad
 A_{03}=-1,\qquad A_{35}=-r,
 \qquad A_{50}=-2a.
\label{eq:four-entries}
\end{equation}
The first closed path has two directed edges, \(1\to2\to1\), and gives
\begin{equation}
\begin{aligned}
 c_{12}
 &=A_{12}A_{21}
   =\left(-\frac{1+s}{2}\right)^2\\
 &=\frac{1+2s+s^2}{4}
   =\frac{6+2s}{4}
   =\frac{3+s}{2},
\end{aligned}
\label{eq:cycle12}
\end{equation}
where \(s^2=5\).  Thus
\begin{equation}
 s=2c_{12}-3\in k(P_{11,j}).
\label{eq:recover-s}
\end{equation}
The second closed path is \(0\to3\to5\to0\):
\begin{equation}
 c_{035}=A_{03}A_{35}A_{50}=(-1)(-r)(-2a)=-2ra.
\label{eq:cycle035}
\end{equation}
Consequently
\begin{equation}
 u:=ra=-\frac12c_{035}\in k(P_{11,j}).
\label{eq:recover-u}
\end{equation}
These equations prove membership in the cyclic-product field.

\section{The radical identity and the complex place}
\label{sec:field-obstruction}

\subsection{Deriving the quadratic relation}

We now prove the one algebraic identity behind all three cases.

\begin{lemma}[Exact radical elimination]
\label{lem:radical-relation}
Let \(s=\sqrt5\), let \(a\) be given by \eqref{eq:a-radical}, and let
\(u=\sqrt2a\).  Then
\begin{equation}
 (s-1)u^2+2(s-3)u+s-7=0.
\label{eq:cycle-relation}
\end{equation}
\end{lemma}

\begin{proof}
Set
\[
 t=\sqrt{2+2s}>0.
\]
By the definitions of \(a\) and \(u\),
\begin{equation}
 u^2=2a^2
 =2\cdot\frac14(7+s+2t)
 =\frac{7+s}{2}+t.
\label{eq:u-square}
\end{equation}
We claim first that
\begin{equation}
 (s-1)u=3-s+t.
\label{eq:linear-radical}
\end{equation}
Both sides are positive: \(s>1\) and \(u>0\) for the left side, while
\(s<3\) and \(t>0\) for the right side.  It is therefore enough to prove
that their squares are equal.  Using \(s^2=5\), \(t^2=2+2s\), and
\eqref{eq:u-square}, expansion gives
\[
\begin{aligned}
 ((s-1)u)^2
 &=(s-1)^2\left(\frac{7+s}{2}+t\right)\\
 &=(6-2s)\left(\frac{7+s}{2}+t\right)\\
 &=16-4s+(6-2s)t,
\end{aligned}
\]
because
\((6-2s)(7+s)/2=(42-8s-2s^2)/2=16-4s\).  On the other hand,
\[
\begin{aligned}
 (3-s+t)^2
 &=(3-s)^2+2(3-s)t+t^2\\
 &=14-6s+2(3-s)t+2+2s\\
 &=16-4s+(6-2s)t.
\end{aligned}
\]
The positive quantities have equal squares, proving
\eqref{eq:linear-radical}.

Rearrange it as
\begin{equation}
 t=(s-1)u+s-3.
\label{eq:t-linear}
\end{equation}
Squaring, substituting \(t^2=2+2s\), and moving everything to the left gives
\[
 ((s-1)u+s-3)^2-(2+2s)=0.
\]
A direct expansion factors the left side as
\[
 (s-1)\bigl((s-1)u^2+2(s-3)u+s-7\bigr).
\]
Because \(s\ne1\), division by \(s-1\) proves the claimed relation.
\end{proof}

\begin{remark}
The expanded proof uses only positivity of the named square roots,
\(s^2=5\), and ring arithmetic.  The Lean theorem
\begin{center}
\code{explicit\_radical\_cycle\_relation}
\end{center}
checks the same identity from the nested square-root definitions.  The
intermediate theorem
\begin{center}
\code{cycle\_relation\_of\_radical\_identities}
\end{center}
checks the final factorization over an arbitrary characteristic-zero field.
\end{remark}

\subsection{Negative discriminant at the conjugate place}

\begin{lemma}[Non-total-reality certificate]
\label{lem:non-real-subfield}
Let \(s=\sqrt5\in\R\) and let \(u\in\C\) satisfy \eqref{eq:cycle-relation}.  Then
\(L=\Q(s,u)\) is not totally real.
\end{lemma}

\begin{proof}
View \eqref{eq:cycle-relation} as \(f(u)=0\), where
\[
 f(X)=(s-1)X^2+2(s-3)X+s-7\in\Q(s)[X].
\]
Its leading coefficient is nonzero because \(s\ne1\).  Its discriminant is
\begin{equation}
\begin{aligned}
 \Delta
 &=\bigl(2(s-3)\bigr)^2-4(s-1)(s-7)\\
 &=4\bigl((s-3)^2-(s-1)(s-7)\bigr)\\
 &=4\bigl((s^2-6s+9)-(s^2-8s+7)\bigr)\\
 &=8(1+s).
\end{aligned}
\label{eq:discriminant}
\end{equation}
Let \(\sigma:\Q(s)\hookrightarrow\R\) be the nonidentity embedding, so
\(\sigma(s)=-s\).  Applying \(\sigma\) to the coefficients gives a real
quadratic \(f^\sigma\) with discriminant
\begin{equation}
 \sigma(\Delta)=8(1-s)<0,
\label{eq:negative-discriminant}
\end{equation}
because \(s=\sqrt5>1\).  By Lemma~\ref{lem:quadratic-test},
\(f^\sigma\) has no real root.

Since \(f\) is a nonzero polynomial with \(f(u)=0\), the extension
\(L=\Q(s,u)\) is a number field.  By
Proposition~\ref{prop:embedding-extension}, \(\sigma\) extends to an
embedding \(\widetilde\sigma:L\hookrightarrow\C\).  Applying this
embedding to \(f(u)=0\) shows that \(\widetilde\sigma(u)\) is a root of
\(f^\sigma\), hence is nonreal.  Therefore \(L\) is not totally real.
\end{proof}

There is a shorter intrinsic formulation.  If a number field \(K\) contains
\(s,u\) with
\[
 \minpoly_\Q(s)=X^2-5,
 \qquad (s-1)u^2+2(s-3)u+s-7=0,
\]
then the root \(-\sqrt5\) of the minimal polynomial is realized as the image
of \(s\) under some embedding \(K\hookrightarrow\C\).  If that embedding
were real, the image of \(u\) would solve the impossible conjugate quadratic.
This is exactly the formulation kernel-checked in
\code{not\_totallyReal\_of\_cycle\_minpoly}.
For the positive radicals used here, \(\Q(s,u)\subset\R\) in its
distinguished embedding.  The conclusion concerns its other embeddings:
at least one has nonreal image.

\section{Proofs of the main theorem and manifold corollary}
\label{sec:main-proof}

\begin{proof}[Proof of Theorem~\ref{thm:main}]
Fix \(j\in\{6,7,8\}\).  Let
\(\Gamma_j=\Gamma(P_{11,j})\) be the group generated by reflections in the
facets of \(P_{11,j}\).  By the Ma--Zheng census statement, \(P_{11,j}\) has
finite volume.  Therefore \(\Gamma_j\) is a lattice.

Let \(K_j=k(\Gamma_j)\) be its adjoint trace field.  By
Theorem~\ref{thm:vinberg-cycles},
\begin{equation}
 K_j=k(P_{11,j}).
\label{eq:vinberg-equality}
\end{equation}
The products in \eqref{eq:cycle12} and \eqref{eq:cycle035} are cyclic
products of \(A_j=2G_j\).  Equations \eqref{eq:recover-s} and
\eqref{eq:recover-u} therefore show that
\[
 s=\sqrt5\in K_j,
 \qquad u=\sqrt2a\in K_j.
\]
All entries of \(G_j\) are algebraic, so their entry field is a number
field, as is its subfield \(K_j\).  Consequently \(L=\Q(s,u)\) is a
number subfield of \(K_j\).

Lemma~\ref{lem:radical-relation} gives the quadratic relation for \(s,u\),
and Lemma~\ref{lem:non-real-subfield} gives a nonreal embedding of \(L\).
By Lemma~\ref{prop:subfield-persists}, that embedding extends to a
nonreal embedding of \(K_j\).  Thus \(K_j\) is not totally real.

If \(\Gamma_j\) were pseudo-arithmetic,
Proposition~\ref{prop:total-real-obstruction} would make \(K_j\) totally
real.  This contradicts the preceding paragraph.  Hence \(\Gamma_j\) is not
pseudo-arithmetic.  The computation was independent of \(j\), so it applies
to all three rows.
\end{proof}

\begin{proof}[Proof of Corollary~\ref{cor:manifolds}]
Choose one \(j\) and write \(\Gamma=\Gamma_{11,j}\).  Let
\(\Gamma^+=\Gamma\cap\Isom^+(\HH^4)\).  The orientation character shows
that \([\Gamma:\Gamma^+]=2\).  As a finite-index subgroup of the finitely
generated reflection group, \(\Gamma^+\) is finitely generated; it is also
linear in characteristic zero.  By Selberg's lemma, it contains a
torsion-free subgroup \(\Lambda\) of finite index.

The action of \(\Lambda\) on \(\HH^4\) is free and orientation preserving,
so
\[
 M=\Lambda\backslash\HH^4
\]
is an orientable hyperbolic \(4\)-manifold.  Since \(\Lambda\) has finite
index in the lattice \(\Gamma\), its quotient has finite volume.  The base
orbifold is noncompact because the census row has at least one ideal vertex;
the finite cover \(M\) is therefore noncompact as explained in
Subsection~\ref{subsec:orbifold-manifold}.

Finite-index invariance of the adjoint trace field gives
\[
 k(\Lambda)=k(\Gamma)=K_j.
\]
Theorem~\ref{thm:main} says that this field is not totally real.  Applying
Proposition~\ref{prop:total-real-obstruction} to \(\Lambda\) shows that
\(\Lambda\), and equivalently the manifold with this fundamental group, is
not pseudo-arithmetic.
\end{proof}

\section{Exact computational certification}
\label{sec:certification}

The preceding proof is hand-checkable after accepting the source rows and
the two external field theorems.  The repository also contains redundant
exact checks.  Their purpose is to detect transcription mistakes and to
make every concrete assertion reproducible.

\subsection{Rank and inertia of all three Gram matrices}

Exact symbolic elimination gives
\[
 \operatorname{rank}(G_j)=5,
 \qquad \operatorname{inertia}(G_j)=(4,1,2)
\]
for \(j=6,7,8\).  Here the entries count positive, negative, and zero
eigenvalues.  The common principal block on facets \(0,1,2,3,4\) has leading
principal minors
\begin{equation}
 1,\quad 1,\quad \frac{5-s}{8},\quad
 \frac{3(5-s)}{32},\quad -\frac{1+s}{32}.
\label{eq:leading-minors}
\end{equation}
They are all nonzero because \(1<s<3\).  The signs are
\(+,+,+,+,-\); the Jacobi signature criterion therefore gives four positive
and one negative direction on this block.  Exact elimination shows that the
full matrix has rank five, so the remaining two directions are null.  The
programs additionally verify that every principal six-by-six minor vanishes.
The nonzero five-by-five minor and the elimination certificate are the rank witnesses; no floating-point eigenvalue is used.

\subsection{Direct finite-volume check for the third polytope}

For \(P_{11,8}\), Lean classifies all \(127\) nonempty principal submatrices:
\(54\) are positive definite, \(12\) are positive semidefinite but not
positive definite, and \(61\) have a negative direction.  It proves that
exactly ten subsets are maximal positive semidefinite, matching the vertex
incidences of the published combinatorial type.  Eight vertices have elliptic rank-four links.
The remaining two are parabolic:
\[
 \{1,3,4,6\}\quad\text{of type }\widetilde C_3,
 \qquad
 \{2,3,4,5,6\}\quad\text{of type }
 \widetilde A_1\sqcup\widetilde C_2.
\]
The second has nullity two and is the Euclidean triangular-prism link of a
nonsimple ideal vertex.  Every elliptic rank-three edge diagram has exactly
two ordinary or ideal endpoints.  Vinberg's finite-volume criterion
\cite{VinbergReflections} then gives finite volume and exactly two cusps.

An independent CoxIter run, pinned to a recorded source commit, returns
\[
 (f_0,f_1,f_2,f_3,f_4)=(10,21,18,7,1),
 \qquad \vol(P_{11,8})=\frac{19\pi^2}{864},
\]
agreeing with the census.  Lean independently proves all five face counts
by closing the vertex--facet incidences and counting nonempty faces,
including the polytope, by strict-chain dimension.  The correspondence
between this incidence poset and geometric faces uses the published
realization; no numerical face count is assumed in the Lean calculation.
CoxIter is described in \cite{Guglielmetti}.
These checks verify the geometric input for one case; the proof of
the family still cites the
Ma--Zheng classification for all three.

\subsection{A complete trace-field check for the third polytope}

The two cycles used in the theorem give only the subfield needed for the
obstruction.  Two further cycles recover \(r=\sqrt2\) and \(b_8\):
\[
 c_{126}=-r(1+s),\qquad c_{03416}=-4b_8.
\]
Since every Gram entry lies in \(\Q(s,r,a,b_8)\), this proves
\[
 k(P_{11,8})=\Q(s,r,a,b_8).
\]
With \(w=2(r+a+b_8)\), exact power-basis reconstruction gives
\(k(P_{11,8})=\Q(w)\), where
\begin{equation}
\begin{aligned}
 p_w(X)={}&X^8-8X^7-252X^6-1392X^5-3524X^4\\
          &-5280X^3-5712X^2-3968X-1136.
\end{aligned}
\label{eq:primitive-polynomial}
\end{equation}
PARI/GP verifies irreducibility, signature
\((r_1,r_2)=(4,2)\), and discriminant
\[
 40960000=2^{16}5^4.
\]
Thus the full field independently exhibits two conjugate pairs of complex
embeddings.  A second primitive-element route produces a different polynomial
for the same explicitly reconstructed field; agreement is not inferred from
degree, signature, and discriminant alone.

\subsection{Artifacts and bounded-memory execution}

The exact programs and durable outputs are:
\begin{itemize}[leftmargin=2em]
\item \path{scripts/p11_6_7_non_total_trace.py};\\
  \path{results/p11_6_7_non_total_trace_certificate.json};
\item \path{scripts/p11_8_non_total_trace.py};\\
  \path{results/p11_8_non_total_trace_certificate.json};
\item \path{scripts/p118_trace_field.py};\\
  \path{results/p118_trace_field_certificate.json};
\item \path{scripts/p118_field_independent.gp}, the independent PARI/GP
  irreducibility, real-root, signature, and discriminant check;
\item \path{scripts/p11_8_finite_volume.py};\\
  \path{results/p11_8_finite_volume_validation.json};
\item the pinned CoxIter input, output, and provenance files under
  \code{results/}.
\end{itemize}
Each expensive symbolic case runs in a fresh process.  On the reference
machine the case-worker resident-set peaks are below \(100\) MB.  The exact
commands, source URLs, versions, and expected digests are in
\code{REPRODUCIBILITY.md}.

The Ma--Zheng source file used by the certifiers has SHA-256
\[
\text{\code{0017931637f017641d4935450cc2a07d070866d9c749443d9e3ede9f3a39f003}}.
\]
Each headline Python certifier checks this digest before parsing the rows.

\section{Lean formalization}
\label{sec:lean}

The directory \code{formal/} is a Lean~4 project pinned to Lean and mathlib
\code{v4.32.0}.  This section states exactly what its successful build means.

\subsection{Unconditional kernel-checked mathematics}

The modules \code{NonTotalReality.lean} and \code{Manuscript.lean}
prove the following.

\begin{enumerate}[leftmargin=2.2em]
\item \code{recover\_s\_from\_twoCycle} and
  \code{recover\_u\_from\_threeCycle} check the rearrangements in
  \eqref{eq:recover-s} and \eqref{eq:recover-u} over every
  characteristic-zero field.
\item \code{explicit\_radical\_cycle\_relation} proves
  \eqref{eq:cycle-relation} from the positive nested-radical definitions
  over \(\R\).
\item \code{minpoly\_eq\_of\_square\_five} proves that \(s^2=5\) implies
  \(\minpoly_\Q(s)=X^2-5\) in any characteristic-zero field.  It uses a
  proof of the irreducibility of \(X^2-5\) over \(\Q\).
\item \code{conjugateQuadratic\_no\_real\_root} proves that the conjugate
  relation has no real solution.  The theorem
  \code{not\_totallyReal\_of\_cycle\_minpoly} then uses mathlib's
  number-field embedding theorem to realize \(s\mapsto-\sqrt5\) and
  deduce non-total reality.  Its square-identity variant is
  \code{not\_totallyReal\_of\_cycle\_square}.
\item A \code{TwoCycleCertificate} records four field elements and four
  equations: the two cycle values, the minimal polynomial of \(s\), and
  the relation for \(u\).  The constructor \code{of\_real\_cycles}
  builds this certificate from a real embedding \(\iota:K\hookrightarrow\R\)
  and the two equalities
  \[
   \iota(c_{12})=(3+\sqrt5)/2,\qquad
   \iota(c_{035})=-2\sqrt2a.
  \]
  It recovers \(s,u\in K\), transports the radical relation by injectivity,
  and proves the minimal-polynomial statement.  The theorem
  \code{not\_totallyReal\_of\_real\_cycles} concludes that \(K\) is not
  totally real.
\end{enumerate}

The module \code{GramAudit.lean} also defines all three displayed matrices
over \(\R\), their entry fields, and their cyclic-product fields.  It proves
the cycle values, algebraicity of the entries, the number-field instances,
and non-total reality of all three cyclic-product fields.  For the third
matrix it also proves equality of the cyclic-product and entry fields.
These results require no geometric predicates.  Identifying these matrices
with the census polytopes, and their cyclic-product fields with the adjoint
trace fields of the reflection groups, remains external.

\subsection{Conditional geometric implications}

The module \code{ExternalContracts.lean} specializes the imported statements
to hyperbolic dimension four.  Its structure \code{CensusInputs} records:
\begin{itemize}[leftmargin=2em]
\item a reflection group \(\Gamma_j\) for each of the three cases, related
  to its actual unit-normal matrix \(G_j\);
\item finite volume and noncompactness of its quotient as separate inputs;
\item Vinberg's equality between the group's adjoint trace field and the
  field generated by closed products in \(2G_j\); and
\item the implication from pseudo-arithmeticity of a finite-volume lattice
  to total reality of its adjoint trace field.
\end{itemize}
The theorem \code{census\_main\_theorem} combines these inputs with the
matrix-derived field proofs to conclude, for every case,
\[
 \text{finite-volume lattice}\ \wedge\
 k(\Gamma_j)=k(G_j)\ \wedge\
 \neg\text{total reality}\ \wedge\
 \neg\text{pseudo-arithmeticity}.
\]
Here \(k(G_j)\) denotes the concrete cyclic-product field.  The trace-field
identification is a literal equality of intermediate fields of \(\R/\Q\).
The theorem requires no separately chosen field or two-cycle certificate.
The group and manifold vocabulary has its cited geometric interpretation;
Lean checks the implications from the stated contracts.

The structure \code{CoverInputs} requires one subgroup--quotient pair
\((\Lambda,M)\) satisfying finite index, torsion-freeness, orientation,
and quotient correspondence simultaneously.  Finite-index invariance is
the equality \(k(\Lambda)=k(\Gamma_j)\).  The field
\code{quotient\_geometry} uses both finite volume and noncompactness of the
base, and supplies the lattice property of \(\Lambda\) as well as the
required manifold geometry.  Thus the pseudo-arithmetic implication applies
to the cover group with its hypotheses intact.  The theorem
\code{census\_manifold\_corollary} proves Corollary~\ref{cor:manifolds}
from these contracts.  The earlier interfaces \code{CoxeterExternalInputs}
and \code{ManifoldExternalInputs} remain reusable abstract implications.

\subsection{Statement-by-statement semantic audit}

The register \path{audit/manuscript_statements.json} maps the mathematical
assertions in every source paragraph, displayed equation, and table to
their Lean declarations or to an explicit external boundary.  Its readable
companion is \path{notes/manuscript_semantic_audit.md}.  Each entry records a
positive example and a negative control: a failed hypothesis, a boundary
case, or a counterexample to a proposed stronger statement.  A valid
theorem has no counterexample satisfying all its hypotheses.  The register
assesses translation separately from proof provenance: external inputs are
acceptable when their specialization, hypotheses, and use match the text.
Reusing a checked library result or a different proof is also acceptable.

The audit modules additionally prove the two embeddings of \(\Q(\sqrt5)\),
embedding extension, the discriminant criterion, the equivalence of total
reality and total positivity for nonzero defining square-root radicands,
the Lorentz reflection identities, and elementary
orientation-index and noncompactness statements.  Exact scripts test all
three matrices and all \(127\) principal subdiagrams of the third, and reject
perturbed cycle values, power-basis coordinates, ranks, and vertex lists.
The five displayed leading-minor identities, the complete principal-subdiagram
classification for \(P_{11,8}\), and its face vector are also proved in Lean.
Further Lean controls reject all \(70\) single-bit vertex--facet incidence
changes and all ten vertex deletions.  They also prove that the incorrect
vectors \((10,21+t,18+t,7,1)\), for every integer \(t>0\), pass the Euler
alternating-sum check.
The contract controls reject an unrelated field, a wrong case--matrix pair,
a reversed implication, and incompatible existential witnesses.  Logical
test models do not instantiate hyperbolic lattices or manifolds; the
geometric interpretation remains the explicit imported input.

\subsection{Other formal modules}

The module \code{NonTotalReality.lean} contains the core field theorem in a
minimal form.  The separate module \code{ConditionalDescent.lean} concerns
the supplementary normalized-Clifford and split-support criterion documented
in the repository notes.  It represents ABHN, local Brauer functoriality,
Hasse--Minkowski, and the pseudo-arithmetic candidate classification as
ordinary hypotheses and checks their downstream implication graph.  That
criterion is not needed for Theorem~\ref{thm:main}; the main examples are
settled by non-total reality before any Clifford invariant is considered.

\subsection{Audit statement}

A release audit searches all project Lean sources for the declarations
\code{axiom} and \code{opaque} and for the proof escapes \code{sorry},
\code{admit}, and \code{native\_decide}.  None occurs in project
declarations or proofs; dependency sources are excluded from this textual
check.  A successful Lean build checks the proof terms against the kernel
and imported library declarations.  The theorem-level axiom audit reports
only Lean's standard \code{propext}, \code{Classical.choice}, and
\code{Quot.sound}.  The geometric conclusions remain conditional on the
hypotheses described above.  The build command is
\[
 \code{cd formal \&\& lake build}.
\]

\section{External-input register}
\label{sec:external-register}

This register records the geometric, field-theoretic, and library inputs.
Each item distinguishes facts proved in the text or in mathlib from those
assumed by the manuscript-facing Lean theorems.

\begin{description}[style=nextline]
\item[\textbf{E1: Census rows and finite volume.}]
Ma--Zheng's classification supplies the three vectors, positive radical
values, cusp counts, volumes, and the assertion that they are finite-volume
hyperbolic Coxeter \(4\)-polytopes \cite{MaZheng}.  Exact scripts verify the
transcription, Gram rank, and inertia; a direct vertex-link proof independently
verifies finite volume for \(P_{11,8}\).  The Lean input
\code{census\_realization} relates each group to its actual matrix, while
\code{finite\_volume} and \code{noncompact} give the two separate geometric
properties.  The census rows themselves remain external data.

\item[\textbf{E2: Vinberg field identification.}]
For a finite-volume Coxeter polytope, the field generated by cyclic products
in \(2G\) is the adjoint trace field of the reflection lattice
\cite{VinbergRings}.  The Lean input \code{vinberg} requires the reflection
group--matrix correspondence and finite volume, and concludes equality
with the concrete cyclic-product field.  The theorem
\code{traceField\_eq\_cycleField} specializes it to each census case.

\item[\textbf{E3: Pseudo-arithmetic fields are totally real.}]
Emery--Mila's definition and Section~1.4 imply that a pseudo-arithmetic
lattice in dimension \(n>3\) has totally real adjoint trace field
\cite{EmeryMila}.  The structure field
\path{pseudoArithmetic_totallyReal} records this implication for
finite-volume lattices in the fixed dimension four; Lean checks its use.

\item[\textbf{E4: Extension and realization of number-field embeddings.}]
The standard extension theorem and mathlib's theorem describing evaluations
of number-field embeddings by roots of a minimal polynomial are used to
realize \(s\mapsto-\sqrt5\).  The extension theorem is recalled in
Proposition~\ref{prop:embedding-extension}; Lean uses mathlib's existing
verified result rather than a project hypothesis.

\item[\textbf{E5: Selberg's lemma.}]
A finitely generated characteristic-zero linear group has a torsion-free
finite-index subgroup \cite{Selberg}.  It is not formalized from first
principles here.  The Lean input \code{cover\_exists} packages the needed
specialization after passing to the orientation subgroup, together with
the quotient existence used in E7.  All properties use the same witnesses.

\item[\textbf{E6: Finite-index trace-field invariance.}]
The adjoint trace field is a commensurability invariant
\cite{VinbergRings}.  The input \path{traceField_invariant} requires
a finite-volume base lattice and a finite-index subgroup, and concludes
equality of their field values.

\item[\textbf{E7: Finite-cover quotient geometry.}]
The standard facts that a torsion-free orientation-preserving subgroup gives
an orientable manifold, finite index preserves finite volume, and a finite
cover of the noncompact quotient remains noncompact are proved at the
elementary level in Subsection~\ref{subsec:orbifold-manifold}.  Lean packages
their geometric realization as \code{quotient\_geometry}, including the
cover's lattice property.  The input \code{quotient\_pseudoArithmetic}
relates the same group and manifold.  These are specialized external
statements; the project does not construct hyperbolic orbifolds.
\end{description}

No item E1--E7 is hidden behind a global Lean axiom.  Replacing an input by a
future library theorem would change the constructor used to build the
corresponding structure, not the checked proof downstream.

\section{Reproducibility and release audit}
\label{sec:reproducibility}

The companion repository~\cite{Repository} contains the exact certificates,
the pinned Lean development, and the reproduction instructions.  The central
conclusion can be checked in three complementary ways.

\paragraph{Level 1: read the certificate.}
Check the four entries in \eqref{eq:four-entries}, multiply the cycles,
follow Lemmas~\ref{lem:radical-relation} and
\ref{lem:non-real-subfield}, and apply Inputs E2 and E3.  No computer algebra
is required.

\paragraph{Level 2: build the formal and written proofs.}
Run
\begin{verbatim}
cd formal
lake build

cd ../arxiv
latexmk -pdf -interaction=nonstopmode -halt-on-error main.tex
\end{verbatim}
The project pins Lean/mathlib, and the TeX source has a fixed date.

\paragraph{Level 3: repeat all exact source checks.}
Follow \code{REPRODUCIBILITY.md} to obtain the pinned Ma--Zheng source and
HCPdm incidence file, verify their SHA-256 values, and run each Python or
PARI/GP job serially.  Reference JSON certificates are committed, so the
mathematical outputs remain inspectable without network access.  Timing and
peak-memory metadata may vary and are not used as mathematical evidence.

The command
\begin{center}
\code{python3 scripts/release\_check.py}
\end{center}
performs the local
release audit: JSON parsing, Python byte-compilation, forbidden-token checks,
documentation consistency, Lean build, exact regressions, and TeX build.  The
packaging script \code{scripts/build\_release.sh} creates a deterministic
archive and checksum under \code{dist/}.  The execution record is split
between \code{PRE\_RELEASE\_CHECKLIST.md} and \code{RELEASE\_NOTES.md}.

\appendix

\section{Dependency graph of the main conclusion}
\label{app:dependency}

For quick auditing, the complete argument is the following chain.
\[
\begin{array}{c}
\text{Ma--Zheng row (E1)}\\[2pt]
\Downarrow\\[-2pt]
\text{four entries of }A=2G\\[2pt]
\Downarrow\\[-2pt]
c_{12}=(3+s)/2,\quad c_{035}=-2u\\[2pt]
\Downarrow\\[-2pt]
s,u\in k(P),\quad \minpoly_\Q(s)=X^2-5\\[2pt]
\Downarrow\\[-2pt]
(s-1)u^2+2(s-3)u+s-7=0\\[2pt]
\Downarrow\\[-2pt]
\text{an embedding extending }s\mapsto-\sqrt5\text{ sends }u\text{ outside }\R\\[2pt]
\Downarrow\\[-2pt]
k(P)\text{ is not totally real}\\[2pt]
\Downarrow\ \text{Vinberg (E2)}\\[-2pt]
k(\Gamma)\text{ is not totally real}\\[2pt]
\Downarrow\ \text{Emery--Mila (E3)}\\[-2pt]
\Gamma\text{ is not pseudo-arithmetic}.
\end{array}
\]
The manifold corollary appends E5--E7.  The radical relation and both
downstream implication chains are checked by Lean; E1--E3 and E5--E7 remain
visible external inputs.

\section{Independent full-field data for comparison}
\label{app:full-field}

The full-field computation is not part of the proof, but it gives a useful
cross-check.  The signature \((4,2)\) of the degree-eight field means that it
has four real embeddings and two conjugate pairs of nonreal embeddings,
since \(r_1+2r_2=8\).  The two-cycle proof detects at least one of those pairs
without constructing the other six embeddings.  The short certificate
therefore proves the required obstruction without a
primitive-element calculation.  The full-field computation gives the
additional degree, signature, and discriminant data.

\section*{Acknowledgments}

The author used OpenAI Codex as a research and coding assistant to help
recover and audit the exact computations, explore the quadratic-form descent
criterion, improve the exposition, and generate and audit the accompanying
Python, PARI/GP, and Lean~4 code.

The author reviewed the resulting mathematics and code and takes
responsibility for the final content.  Codex is acknowledged as a tool, not
listed as an author.

\bibliographystyle{plain}
\begingroup
\raggedright
\bibliography{references}
\endgroup

\end{document}